\documentclass[11pt]{article}
\usepackage{amsmath}
\usepackage{amsfonts}
\usepackage{amssymb}
\usepackage{latexsym}
\usepackage[mathscr]{eucal}
\usepackage{graphicx}
\usepackage{epsf}
\usepackage{epsfig}
\usepackage[english]{babel}
\usepackage[applemac] {inputenc}
\usepackage[T1]{fontenc}

\usepackage[usenames]{xcolor}

\def\rnew{\color{red} }

\usepackage{fullpage}
\def\RR{\rm \hbox{I\kern-.2em\hbox{R}}}
\def\NN{\rm \hbox{I\kern-.2em\hbox{N}}}
\def\ZZ{\rm {{\rm Z}\kern-.28em{\rm Z}}}
\def\CC{\rm \hbox{C\kern -.5em {\raise .32ex \hbox{$\scriptscriptstyle
|$}}\kern
-.22em{\raise .6ex \hbox{$\scriptscriptstyle |$}}\kern .4em}}

\def\<{\langle}
\def\>{\rangle}
\def\t{\tilde}

\def\e{\varepsilon}

\def\nl{\newline}
\def\o{\overline}
\def\wt{\widetilde}

\def\cT{{\cal T}}

\def\cL{{\cal L}}
\def\cA{{\cal A}}

\def\cK{{\cal K}}

\def\cM{{\cal M}}
\def\cF{{\cal F}}
\def\cP{{\cal P}}

\def\cP{{\cal P}}
\def\cO{{\cal O}}
\def\R{\mathbb{R}}
\def\N{\mathbb{N}}

\def\Chi{\raise .3ex
\hbox{\large $\chi$}} 
\def\lsima{\hbox{\kern -.6em\raisebox{-1ex}{$~\stackrel{\textstyle<}{\sim}~$}}\kern -.4em}
\def\lsim{\hbox{\kern -.2em\raisebox{-1ex}{$~\stackrel{\textstyle<}{\sim}~$}}\kern -.2em}
\def\gsim{\hbox{\kern -.2em\raisebox{-1ex}{$~\stackrel{\textstyle>}{\sim}~$}}\kern -.2em}
\def\[{\Bigl [}
\def\]{\Bigr ]}
\def\({\Bigl (}
\def\){\Bigr )}
\def\[{\Bigl [}
\def\]{\Bigr ]}
\def\({\Bigl (}
\def\){\Bigr )}

\def\span{\mathop{span}}

\newcommand{\be}{\begin{equation}}
\newcommand{\ee}{\end{equation}}
\newcommand{\bea}{$$ \begin{array}{lll}}
\newcommand{\eea}{\end{array} $$}
\newcommand{\bi}{\begin{itemize}}
\newcommand{\ei}{\end{itemize}}
\newcommand{\iref}[1]{(\ref{#1})}
\newtheorem{theorem}{Theorem}
\newtheorem{remark}{Remark}
\newtheorem{lemma}{Lemma}

\newtheorem{corollary}{Corollary}

  \newcommand{\wnote}[1]{\marginpar{\tiny \textcolor{red}{WD: #1}}}
  
\def\P{\mathbb P}
\def\N{\mathbb N}

\def\argmax{\mathop{\rm argmax}}

\def\dist{\mathop{\rm dist}}
\def\span{\mathop{\rm span}}

\def\anew{\color{black}}

\newcommand{\surr}{\bar{e}}

\begin{document}
\bibliographystyle{plain}
\title{
  Reduced Basis Greedy Selection Using Random Training Sets
  }
\author{ 
 Albert Cohen, Wolfgang Dahmen,
 Ronald DeVore, and James Nichols
\thanks{%
   This research was supported by the ONR Contracts
  N00014-11-1-0712,  N00014-12-1-0561,  {and} N00014-15-1-2181; 
  the  NSF Grants  DMS 1222715,    DMS 1720297,
       DMS 1222390,  and DMS 1521067;
  the Institut Universitaire de France; The European Research Council under grant ERC AdG 338977 BREAD; 
  and by the SmartState and Williams-Hedberg Foundation.
  }  
  }
 \maketitle
\date{}
 

\begin{abstract} 
Reduced bases have been introduced for the approximation of
parametrized PDEs in applications where many online queries are
required. Their numerical efficiency for such problems has been 
theoretically confirmed in \cite{BCDDPW,DPW}, where it is shown that 
the reduced basis space  $V_n$ of dimension $n$, constructed by a certain greedy strategy, has  approximation error  similar to that of the optimal space associated to the Kolmogorov $n$-width of the 
solution manifold. The greedy construction of the reduced basis space is performed
in an offline stage which requires at each step a maximization of the current error over 
the parameter space. For the purpose of numerical computation, this maximization
is performed over a finite {\em training set} obtained through a discretization of  the parameter domain.
To guarantee a final approximation error $\e$ for the space generated by the greedy algorithm requires 
in principle that the snapshots associated to this training set constitute an approximation net 
for the solution manifold with accuracy of order $\e$.   Hence, the size of the training set 
is the $\e$ covering number for $\cM$ and this covering number typically behaves like
$\exp(C\e^{-1/s})$ for some $C>0$ when the solution manifold has $n$-width decay $O(n^{-s})$.  
Thus, the shear size of the training set prohibits implementation of the algorithm when $\e$ is  small. 
The main result of this paper shows that, if one is willing to accept results which hold with high probability, rather than with  certainty, then for
a large class of relevant problems one may  replace the fine discretization by  a random training set of  size polynomial in $\e^{-1}$.   Our proof of this fact is established by using inverse
inequalities for polynomials in high dimensions.  
\end{abstract}

%

\section{Introduction}

Complex systems are frequently described by parametrized PDEs that take the general form
\be
\cP(u,y)=0.
\ee
Here $y=(y_j)_{j=1,\dots,d}$ is a vector of parameters ranging
over  some domain $Y\subset \R^d$ and $u=u(y)$ is the corresponding solution
which is assumed to be uniquely defined  in some Hilbert space $V$ 
for every $y\in Y$. {\anew Throughout this paper}, we denote by $\|\cdot\|$ and $\<\cdot,\cdot\>$ the norm and inner product of $V$.  In what follows, we assume that the parameters have been rescaled so that $Y=[-1,1]^d$.  Here $d$ is typically large and in some cases $d=\infty$.   We seek results that are immune to the size of $d$, i.e., are dimension independent.

 Various {\it reduced modeling} approaches have been developed for the purpose 
of efficiently approximating the solution $u(y)$ in the context
of applications where the solution  map
\be
y\mapsto u(y),
\ee
needs to be queried for a large number of parameter values $y\in Y$.
This need occurs for example in optimal design or inverse problems
where such parameters need to be optimized. 
The strategy consists in first constructing in some {\it offline stage}
a linear space $V_n$ of hopefully low dimension $n$,
that provides  a reduced map
\be
y\mapsto u_n(y) \in V_n, 
\ee
that  approximates the solution map to the required target accuracy $\e$ for all queries of $u(y)$.  The reduced map is then implemented in the {\it online stage} with greatly reduced computational cost, typically polynomial in $n$.
 
As opposed to standard approximation spaces such as finite elements, the spaces
$V_n$ are  specifically designed   to approximate the image of $u$, i.e., the elements in  the parametrized family
\be
\cM:=\{u(y)\; : \; y\in Y\},
\ee
which is called the {\it solution manifold}.   The optimal  $n$-dimensional linear approximation space
for this manifold is the one that achieves the minimum in the definition of the 
Kolmogorov $n$-width
\be
d_n=d_n(\cM)_V:=\inf_{\dim(E)=n} \sup_{u\in \cM} \|u-P_E u\|.
\ee
This optimal space is computationally out of reach and the above quantity should be viewed as a
benchmark for more practical methods.
Here $P_E$ denotes again  the $V$-orthogonal projection onto $E$ for any subspace $E$ of $V$.

One approach for constructing a reduced space $V_n$, which comes with
substantial theoretical footing and will be instrumental in our discussion,
consists in proving that the solution map is analytic in the parameters $y$ and has a Taylor expansion%
\be
\label{te}
u(y)=\sum_{\nu\in\cF}t_\nu y^\nu,
\ee
where $y^\nu:=\prod_{j=1}^d y_j^{\nu_j}$
and $\cF:=\{\nu\in \N^d\}$. In the case of countably many parameters, 
$\cF$ is the set of finitely supported sequences $\nu=(\nu_j)_{j\geq 1}$ with $\nu_j\in\N$.
One then proves that the sequence $(t_\nu)_{\nu\in\cF}$ of Taylor coefficients has some decay property.  Two prototypical examples of results on decay are given in \cite {CDS,BCM} for the elliptic equation
\be
\label{ellip}
{\rm div}(a\nabla u)=f,
\ee
where the diffusion coefficient function has the parametrized form
\be
a=a(y)=\o a+\sum_{j\geq 1}y_j\psi_j,
\label{affine}
\ee
for some given functions $\o a$ and $(\psi_j)_{j\geq 1}$. 
These results {\anew (see e.g. Theorem 1.1 and 1.2 in \cite{BCM})} show that, under mild decay or summability conditions on the
functions $\psi_j$, one has that the sequence $(\|t_\nu\|_V)_{\nu\in \cF}$ is in $\ell^p$ for certain 
$p<1$ with a bound on the $\ell^p$ {\anew quasi-norm}.  It then follows that for each $n$, there is a  set $\Lambda_n\subset \cF$ with 
$\#(\Lambda_n)=n$ such that 
\be
\label{te2}
\sup_{y\in Y} \|u(y)-\sum_{\nu\in\Lambda_n}t_\nu y^\nu\|\le Cn^{-r},\quad r:=1/p-1.
\ee
Similar results have been obtained for more general models of linear and nonlinear PDEs, see 
in particular \cite{CCS,CD}{\anew, using orthogonal expansions into tensorized Legendre polynomials
in place of Taylor series.}

Therefore the space $V_n:={\rm span}\{ t_\nu \; : \; \nu\in\Lambda_n\}$ has dimension $n$ with an a priori bound $Cn^{-r}$ on its approximation error for all members of the solution manifold $\cM$. One choice of $\Lambda_n$
giving \iref{te2} is the set of indices corresponding to the $n$ largest $\|t_\nu\|_V$. Further analysis \cite{CM} shows that
the same convergence estimate can be obtained imposing in addition that the sets
$\Lambda_n$ are {\it downward closed} (or {\it lower sets}), i.e. having the property 
\be
\nu\in\Lambda_n \quad{\rm and} \quad \mu\le \nu \implies \mu\in \Lambda_n,
\ee
where $\mu\le \mu$ is to be understood
componentwise. We stress that the rate of decay $n^{-r}$ in the bound \iref{te2} may be suboptimal compared to 
the actual rate of decay of the $n$-width $d_n(\cM)$.

The present paper is concerned with another prominent reduced modeling strategy known
as the {\it Reduced Basis Method} (RBM). In this approach \cite{MPT,MPT1,RHP, S},
particular snapshots
\be
u^k=u(y^k), \quad k=1,2,\dots,
\ee
are selected in the solution manifold and the space $V_n$ is defined by
\be
V_n:= {\rm span} \{u^1,\dots,u^n\}.
\ee
A certain {\it greedy procedure} has been {\anew first proposed in \cite{PPRV} and 
analyzed} in \cite{BMPPT} for selecting these snapshots.
It was shown in \cite{BCDDPW,DPW} that  the approximation error
\be
\sigma_n=\sigma_n(\cM)_V:=\sup_{u\in \cM} \|u-P_{V_n} u\|,
\ee
provided by the resulting spaces has the {same rate of decay (polynomial or exponential) as that of $d_n$. 
In this sense, the method leads to reduced models with optimal performance, in contrast to
sparse polynomial expansions.}

 In its simplest (and idealized) form,  the greedy procedure can be described as follows:
at the initial
step, one sets $V_0=\{0\}$, and given that $V_n$ has been produced after $n$ steps, 
one selects the new snapshot by
\be
u^{n+1}:= \argmax_{u\in \cM} \|u-P_{V_n} u\|. 
\label{greedy}
\ee
Each greedy step thus amounts to maximizing  
\be
e_n(y):=\|u(y)-P_{V_n} u(y)\|
\label{eny}
\ee
over the parameter 
domain $Y$. While in this precise form the scheme cannot be realized in practice an important modification of this greedy selection, known as {\it the weak greedy algorithm},  allows the selection to be done in a a practically feasible  manner while retaining the same performance guarantees, see \S 
\ref{sec:performance}
below.  

 The optimization in the greedy algorithm  is typically performed by replacing $Y$  at each step by a 
discrete {\it training set} $\t Y$.  In order  to retain the performance guarantees of the greedy algorithm,  
this discretization should in principle be chosen fine enough so that the solution
map $y\mapsto u(y)$ is resolved up to the target accuracy $\e>0$, that is,
the discrete set 
\be
\wt \cM=\{u(y)\; : \; y\in \t Y\},
\label{disctrain}
\ee
is an $\e$-approximation net for $\cM$. 

Although performed in the offline stage, this discretization becomes computationally
problematic when the parametric dimension $d$ is either large
or infinite, due to the prohibitive size of this net as $\e\to 0$.   For example, in the typical case
when the Kolmogorov width of $\cM$ decays like $O(n^{-s})$ for some $s>0$, 
we can invoke Carl's inequality \cite{Pi} to obtain a sharp bound
$e^{c{\e^{-1/s}}}$ for the cardinality of $\t \cM$ and $\t Y$. 
{\anew This exponential growth drastically limits the possibility of using
$\e$-approximation nets in practical applications when the number of involved
parameters become large}.

There is a preference toward the use of the greedy constructions over the Taylor expansion constructions because they guarantee error decay 
comparable to the decay of Kolmogorov
widths while the Taylor polynomial constructions do not provide any such guarantee.   Therefore, it is of interest to understand whether the apparent impediment of requiring such a fine discretization of
the solution manifold can somehow be avoided or significantly mitigated. The main result of this paper is to  prove that this is indeed the case provided that one is willing to accept error guarantees that hold with {\em high probability}
rather than with certainty.   Our main result shows that a target accuracy $\e$ can generally be met with high probability by searching over a randomly discretized set $\t Y$ 
whose size grows only polynomially in $\e^{-1}$ rather than exponentially, {\anew 
in contrast to $\e$-approximation nets}.

The   paper is organized as follows. In  \S \ref{sec:performance} , we elaborate on the  weak form of the
greedy algorithm, which is used in numerical computation, 
and recall some known facts on its performance and complexity.
In \S \ref{sec:surrogate}, we use properties of downward closed 
polynomial approximation to show how a random sampling $\t Y$ 
provides  an approximate solution of 
optimization problem engaged at each step of the greedy algorithm.   
In \S \ref{sec:main}, we formulate our modification of the greedy algorithm 
based on such random selection and then analyze its performance. {\anew Finally
we illustrate in \S \ref{sec:test} the validity of the randomized approach by some numerical tests performed
in parametric dimensions up to $d=16$ for which the size of $\e$-approximation nets become
computationally prohibitive.}

\section{Performance and complexity of reduced basis greedy algorithms}
\label{sec:performance}

The greedy selection process described in the introduction is not practically feasible, due to 
at least three obstructions:
\begin{enumerate}
\item
Given a parameter value $y$, the
snapshot $ u(y)$, in particular, the generators of the reduced spaces cannot be
exactly computed. 
\item
For a given $y\in Y$, the quantity $e_n(y)$
to be maximized cannot be exactly evaluated.
\item
The map $y\mapsto e_n(y)$ 
is non-convex/non-concave
and therefore difficult to maximize, even if it could be exactly evaluated.
\end{enumerate}

The first obstruction can be handled when a numerical solver is available for computing 
an approximation
\be
y\mapsto u_h(y),
\ee
of $u(y)$ to any prescribed accuracy $\e_h>0$, that is, such that
\be
\sup_{y\in Y} \|u(y)-u_h(y)\|\leq \e_h.
\label{discerr}
\ee
Here $h>0$ is a space discretization parameter: typically, $u_h$ belongs to a finite element
space $V_h$ of meshsize $h$ and (possibly very large) dimension $n_h$. The selected reduced basis functions 
are now given by $u_i=u_h(y^i)\in V_h$ and therefore the reduced
basis space $V_n$ is a subspace of $V_h$. Whenever the $n$-widths $d_n$ decay much faster than the 
approximation order provided by $V_h$  the reduced space 
$V_n$ has typically much smaller dimension than $V_h$, that is
\be
n<< n_h.
\ee
This yields substantial computational savings when using the reduced basis 
discretization in the online stage.

Note that this numerical solver allows us in principle to also handle the
second obstruction: we could now perform the greedy algorithm 
by maximizing at each step the quantity
\be
e_{n,h}(y):=\|u_h(y)-P_{V_n} u_h(y)\|,
\label{enh}
\ee
which, in contrast to $e_n(y)$, can be exactly computed and satisfies
\be
|e_n(y)-e_{n,h}(y)| \leq \e_h, \quad y\in Y.
\label{enenh}
\ee 
In other words, the greedy algorithm is applied on the approximate
solution manifold
\be
\cM_h:=\{u_h(y)\; : \; y\in Y\}.
\ee
While the quantity $e_{n,h}(y)$ can in principle be computed exactly, 
the complexity of this computation depends, at least in a linear 
manner, on the dimension $n_h={\rm dim}(V_h)$, which is typically much higher than $n$. 
Substantial computational saving may still be obtained when 
maximizing instead an a-posteriori estimator $\o e_{n,h}(y)$ of this quantity
that satisfies
\be
\alpha \o e_{n,h}(y) \leq e_{n,h}(y) \leq \beta \o e_{n,h}(y),\quad y \in Y.
\label{frameh}
\ee
{\anew The computation of $\o e_{n,h}(y)$ for a given $y\in Y$ is based 
in particular on replacing the orthogonal projection $P_{V_n} u_h(y)$
by a Galerkin projection. In turn, it does {\em  not} require the computation of $u_h(y)$
and entails a computational cost $c(n)$ depending on the small dimension $n$, typically
in a polynomial way, rather than on the large dimension $n_h$.}
We refer to \cite{CD} for the derivation of a residual-based estimator $\o e_{n,h}(y)$ having 
these properties in the case of elliptic 
PDEs with affine parameter dependence. 

Maximizing the a-posteriori estimator $\o e_{n,h}(y)$
amounts in applying to $\cM_h$ a so called {\it weak-greedy} algorithm,
where $u^{n+1}$ now satisfies
\be
\|u^{n+1}-P_{V_n} u^{n+1}\|\geq \gamma \|u-P_{V_n} u\|,\quad u\in \cM_h,
\ee
with parameter $\gamma:=\alpha/\beta\in ]0,1[$. 

For such an algorithm, it was proved in \cite{BCDDPW,DPW} that any polynomial or exponential rate of decay
achieved by the Kolmogorov $n$-width $d_n$ is retained by the error performance $\sigma_n$ for this algorithm. More precisely, the following holds, for any compact set $\cK$
in a Hilbert space $V$, see \cite{CD}.

\begin{theorem}
\label{theo1}
Let $d_n=d_n(\cK)_V$ be the $n$-widths of the solution manifold.
Consider the weak greedy algorithm with threshold parameter $\gamma$. 
  For any $C_0>0$
and $s>0$,  we have
\be
\label{Cs}
d_{n}\le C_0(\max\{1,n\})^{-s}, \quad n\ge 0\quad \implies \quad \sigma_n\le C_s \gamma^{-2}(\max\{1,n\})^{-s}, \quad n\ge 0,
\ee 
where $C_s:=  2^{4s+1}C_0$.   For any $c_0,C_0>0$ and  $s>0$,  we have
\be 
d_{n}\le C_0e^{-c_0n^{s}}, \quad n\ge 0\quad \implies \quad \sigma_n\le C_s \gamma^{-1}e^{-c_1n^{s}}, \quad n\ge 0,
\ee 
where $c_1=\frac {c_0} 2 3^{-s}$ and $C_1:=C_0\max\{\sqrt{2},e^{c_1}\}$. 
\end{theorem}

\begin{remark}
\label{remrange}
If the same rates of $d_n$ in the above theorem are 
only assumed within a limited range $0\leq n\leq n^*$, then the same
decay rates of $\sigma_n$ are achieved for the same range $0\leq n\leq n^*$,
up to some minor changes in the expressions of the constants $c_1$ and $C_1$,
independently of $n^*$.
\end{remark}

{\anew
\begin{remark}
\label{necessaryrate}
The reduced basis algorithm aims to construct a $n$-dimensional linear
space that is taylored for the approximation of {\em all} solutions that constitute the solution manifold $\cM$.
Therefore, its performance $\sigma_n$ is always bounded by below by the $n$-width $d_n=d_n(\cM)$.
Assuming some rate of decay on $d_n$, at least polynomial, is thus strictly
necessary is we want the reduced basis method to converge at such a rate.
As discussed in the introduction, in the high dimensional parametric context,
such rate can be rigourously established for certain instance of linear elliptic PDE's 
with affine parametrization of the diffusion function. 
A more general approach applicable to nonlinear PDE's and non-affine parametrizations
for proving such rates is given in \cite{CD}. On the other hand, it is known that $d_n(\cM)$ has poor decay
for certain categories of 
parametrized problems. This  includes in particular transport dominated problems with sharp transition
locus that varies with value of the parameter $y$. Such problems
are thus intrisically not well taylored to reduced basis methods.
\end{remark}
}

The additional perturbations due to the numerical solver and the 
a-posteriori error indicator can thus be incorporated in the analysis
of the reduced basis algorithm. If $\e>0$ is our final target accuracy, 
we set the space discretization parameter $h$ so that $\e_h=\frac \e 2$ {\anew where $\e_h$ is the
space discretization error bound in \iref{discerr}}. We then apply the greedy selection 
on $\cM_h$ based on maximizing $\o e_{n,h}(y)$
until we are ensured that $e_{n,h}(y)\leq \frac \e 2$ for all $y\in Y$. The target accuracy
$e_n(y)\leq \e$ is thus met for all $y\in Y$ in view of \iref{enenh}.

Note that a decay rate $d_{n}(\cM)\leq \gamma(n)$ for some decreasing sequence $\gamma(n)$ 
implies a comparable rate $d_{n}(\cM_h)\leq 2\gamma(n)$, for the range $n\leq n^*$ where 
$n^*$ is the largest value of $n$ such that $\gamma(n)\geq \e_h$. Therefore, using Remark \ref{remrange}
in conjunction with Theorem \ref{theo1} applied to $\cM_h$, we obtain an estimate
on the number of greedy steps $n(\e)$ that are necessary to reach the target accuracy $\e$.

\begin{corollary}
\label{cordnh}
Let $d_n=d_n(\cM)_V$ be the $n$-width of the solution manifold.
For any $s>0$ and $C_0>0$,
\be
d_{n}\le C_0(\max\{1,n\})^{-s},\quad n\geq 0 \quad \implies n(\e) \leq C_1\e^{-\frac 1 s}, \quad \e>0,
\ee
where $C_1$ depends on $C_0$, $s$ and $\gamma$.  For any $s>0$ and $c_0,C_0>0$, 
\be 
d_{n}\le C_0e^{-c_0n^{s}}, \quad n\ge 0\quad \implies \quad n(\e) \leq C_1\max\{\log(c_1\e)^{1/s},0\}, \quad \e>0,
\ee 
where $C_1$ and $c_1$ depend on $C_0$, $c_0$, $s$ and $\gamma$. 
\end{corollary}

The difficulty in item 3 is the most problematic one, in particular when the parametric 
variable $y$ is high-dimensional,  and is the main motivation for the present work.
Since the quantities   $e_n(y)$, $\o e_n(y)$, $e_{n,h}(y)$ and $\o e_{n,h}(y)$ 
may have many local maxima, continuous optimization
techniques are not appropriate. A typically employed strategy is therefore  to replace the continuous optimization
over $Y$ by its discrete optimization over a training set $\wt Y\subset Y$
of finite size. This amounts to applying the greedy or weak-greedy algorithm to
the discretized manifold  $\wt \cM$ {\anew defined by \iref{disctrain}}
or, more practically, to its approximated version
\be
\wt \cM_h:= \{u_h(y) \; : \; y\in \wt Y\}.
\ee

On a first intuition, the discretization should be sufficiently fine so that the manifold $\cM_h$
is resolved with accuracy of the same order as the target accuracy $\e$.  Recall that
 if $K$ is a compact set in some normed space, a finite set $S$ is called a $\delta$-net
of $K$ if 
\be
K\subset {\displaystyle  \bigcup_{v\in S}}B(v,\delta),
\ee
that is, any $u\in K$ is at distance at most $\delta$ from some $v\in S$. 

 The perturbation of the greedy algorithm due to this discretization can be  accounted for jointly with the
previously discussed perturbation, namely  finite element approximation and a-posteriori
error estimation.  Assuming for example that $\wt \cM_h$ is a $\e/3$-net of $\cM_h$,
we set the space discretization parameter $h$ so that
$\e_h=\frac \e 3$. We then apply the greedy selection 
on $\wt \cM_h$ based on maximizing $\surr_{n,h}(y)$   over $\wt Y$
until we are ensured that $e_{n,h}(y)\leq \frac \e 3$ for all $y\in \wt Y$.
  By the covering property, we have $e_{n,h}(y)\leq 2\e/3$ for all $y\in Y$,
and therefore the target accuracy $e_n(y)\leq \e$ is met for all $y\in Y$.
In addition, since $d_n(\wt \cM_h)_V\leq d_n(\cM_h)_V$, 
the statement of Corollary \ref{cordnh} remains unchanged for this discretized algorithm.

The main problem with this approach is that the size of an $\e$-net
of $\cM$ or $\cM_h$ becomes extremely large, especially in high parametric dimension. 

A first natural strategy to generate 
such an $\e$-net  would be to apply the solution map $y\mapsto u(y)$ 
to an $\e$-net for $Y$ in a suitable norm, relying on a stability estimate for this map.
For example, in the simple case of the elliptic PDE \iref{ellip}
with parametrized coefficients \iref{affine}, one can easily establish a stability estimate of the form
{\anew
\be
\|u(y)- u(\t y)\|\leq C\|y-\t y\|_{\ell^\infty}, \quad y,\t y\in Y,
\ee
under the minimal uniform ellipticity assumption $\sum_{j\geq 1} |\psi_j|\leq \min \o a-\delta$ for some $\delta>0$, with $C$ depending on $\min \o a$ and $\delta$.}
Thus an $\e$-net of $\cM$ or $\cM_h$ in the $V$ norm 
is induced by a $C^{-1}\e$-net $\wt Y$ of $Y$ in the $\ell^\infty$ norm.
However, the size of such a net scales like 
\be
\label{coverscale0}
\#(\wt Y)\sim \e^{-cd}.
\ee
with the parametric dimension $d$, therefore suffering
from the curse of dimensionality. In the  case $d=\infty$, one would have to truncate
the parametric expansion \iref{affine} for a given target accuracy $\e$, {\anew resulting 
into an active parametric dimension $d(\e)<\infty$.
Assuming a polynomially decaying error   
$\|\sum_{j>k}|\psi_j|\|_{L^\infty}\lsim k^{-b}$, the growth of $d(\e)$ 
as $\e\to 0$ is in $\cO(\e^{-1/b})$} resulting in a 
training set of size scaling like 
\be
\label{coverscale1}
\#(\wt Y)\sim \e^{-c\e^{-1/b}},
\ee
which is extremely prohibitive.

One sharper way to obtain an estimate independent of the parametric dimension 
is to use a fundamental result that relates covering and widths. We define
the entropy number $\e_n:=\e_n(\cM)_V$ as the smallest value of $\e>0$ such that 
there exists a covering of $\cM$ by $2^n$ balls of radius $\e$. Then, Carl's inequality \cite{Pi}
states that for any $s>0$,
\be
(n+1)^s \e_n \leq C_s\sup_{k=0,\dots,n} (k+1)^s d_k,
\ee
where $d_k=d_k(\cM)_V$ and $C_s$ is a fixed constant.
This inequality shows that, in the case of polynomial decay $n^{-s}$ of 
the $n$-widths, there exists an $\e$-net $\wt \cM$ associated with a training set of size
\be
\label{coverscale}
\#(\wt Y)\sim 2^{c\e^{-1/s}}.
\ee
While this estimate is more favorable than \iref{coverscale1}, it is still extremely prohibitive.
Moreover, the construction of such an $\e$-net, as in the proof of Carl's inequality, necessitates
the knowledge of the approximation spaces $V_n$ that perform with the $n$-width accuracy
$n^{-s}$ which is precisely the objective of the greedy algorithm.

The computational cost at each step $n$ of the offline stage is determined by the product 
between $\#(\wt Y)$ and the cost $c(n)$ of evaluating the error bound $\o e_{n,h}(y)$ for an individual $y\in \wt Y$. 
Therefore, the prohibitive number of error bound evaluations 
is the limiting factor in practice and poses the main obstruction 
to the feasibility of certified reduced basis methods in the regime of
polynomially decaying $n$-widths, and hence in particular, in the context of high parametric dimension.

  In what follows, we show that this obstruction can be circumvented
by not searching for an $\e$-net of $\cM$ but rather defining   $\wt Y$
by random sampling of $Y$. This approach allows us 
to significantly reduce the size of training sets used in greedy algorithms
while still obtaining reduced bases with the same guarantee of performance at least with high probability.   
In order to keep our arguments and notation as simple and clear as possible, we do not consider the issue 
of space discretization and error estimation, assuming that we have access to $e_n(y)$ for each individual $y\in Y$.  
As just described, a corresponding finer analysis  
can incorporate the perturbation of using instead $\o e_n(y)$, $e_{n,h}(y)$ or $\o e_{n,h}(y)$,
with the same resulting overall performance.

\section{Polynomial approximation}
\label{sec:surrogate}

Let $V_n:=\span\{u^1,\dots,u^n\}$ be the reduced basis space at the $n$-th step of the weak greedy algorithm. The next step of the greedy algorithm is  to search over $Y$ to find
a point $y\in Y$ where 
\be
\label{errorn}
e_n(y):=\|u(y)-P_{V_n}u(y)\|_V
\ee
(in practice $\surr_{n,h}(y)$) is large, hopefully close to its maximum over $Y$.   In this section we show that random sampling gives a discrete set  $\wt Y$, of moderate size,  on which the maximum of $e_n(y)$ can be compared with the maximum
of $e_n(y)$ over all of $Y$ with high probability.  To obtain a result of this type we use approximation by polynomials.

Recall that $\Lambda\subset \cF$ is said to be a {\it downward closed set} if whenever $\nu\in\Lambda$ and $\mu\le \nu$, then $\mu\in \Lambda$, where $\mu\le \mu$ is to be understood
componentwise.  To such a set $\Lambda$, we associate
the multivariate polynomial space
\be
\P_\Lambda:={\rm span}\{y\mapsto y^\nu:=\prod_{j\geq 1}y_j^{\nu_j}\; : \; \nu\in\Lambda\}.
\ee
We define 
\be
\cP_\Lambda:=V\otimes \P_\Lambda,
\ee
the space of $V$-valued polynomials
spanned by the same monomials.   
Thus, any polynomial $P$ in $\cP_\Lambda$ takes 
the form $P(y)=\sum_{\nu\in\Lambda}a_\nu y^\nu$ where the $a_\nu$ are in $V$.   For any $m\ge 0$, we  let 
\be
\label{sigma}
\Sigma_m:=\bigcup_{\#(\Lambda)=m}\cP_\Lambda,
\ee 
where the union is over all downward closed sets of size $m$. {\anew Note that the union of all downward
closed sets of cardinality $m$ is the so-called hyperbolic cross $H_{m,d}$ that consists of 
all $\nu$ such that $\prod_{j=1}^d (1+\nu_j) \leq m$.}

Given a function $ v$ in $L^\infty(Y,V)$, we consider its approximation in $L^\infty(Y,V)$ by the elements of $\Sigma_m$   and the error
\be
\label{esigma}
\delta_m(v):= \inf_{P\in\Sigma_m} \sup_{y\in Y}\|v(y)-P(y)\|_V.
\ee
For $r>0$, a function $ v$ in $L^\infty(Y,V)$ is said to be in the approximation class $\cA^r = \cA^r((\Sigma_m)_{m\geq 1})$ if 
\be
\label{ac}
\delta_m(v)\le Cm^{-r}, \quad m\ge 1,
\ee
and the smallest such $C$ defines  a quasi-seminorm $|v|_{\cA^r}$ for this class which is
linear subspace of {\anew $L^\infty(Y,V)$}. A quasi-norm for this space is defined by
\be
\|v\|_{\cA^r}:=\max \{ \|v\|_{L^\infty(Y,V)},|v|_{\cA^r}\}.
\ee
Several foundational results in parametric PDEs prove  that the solution map $y\mapsto u(y)$ belongs
to classes $\cA^r$, as already mentioned in our introduction.
 An important observation to us is that whenever $u\in \cA^r$ and $V_n$ is a finite dimensional subspace of
$V$ then both $P_{V_n}u$ and $u-P_{V_n}u$ are also in $\cA^r$.  For example, for any downward closed set $\Lambda\subset \cF$ and an approximation $P(y)=\sum_{\nu\in\Lambda}a_\nu y^\nu$ to $u(y)$,
the polynomial $Q(y):=\sum_{\nu\in\Lambda} P_{V_n}a_\nu y^\nu$ is in $\cP_\Lambda$ and
\be
\label{observe}
\|P_{V_n}u(y)-Q(y)\|_V=\|P_{V_n}(u(y)-P(y))\|_V \le \|u(y)-P(y)\|_V.
\ee
 It follows that  $\delta_m(P_{V_n}u)\le \delta_m(u)$ for all $m\geq 1$.  The same holds for $u-P_{V_n}u=P_{V_n^\perp}u$.
 From this, one derives that 
\be
\label{boundnorms}
\|P_{V_n}u\|_{\cA^r}\le \|u\|_{\cA^r}\quad {\rm and}\quad \|u-P_{V_n}u\|_{\cA^r}\le \|u\|_{\cA^r}.
\ee 

The next result shows that when a function belongs to the class $\cA^r$, its maximum
over a random set of point $\wt Y$ is above a fixed fraction of its maximum over $Y$
with some controlled probability.

\begin{lemma}
\label{slemma}
Let $r>1$ and suppose that $v\in \cA^r$ with $\|v\|_{\cA^r}\le M_0$ and $\|v\|_{L^\infty(Y,V)}=M$.   Let  $m$ be
an integer such that       
\be
\label{mcond}
 4 M_0m^{-r+1}< M.
  \ee
If $\wt Y$ is any finite set of $N$ points drawn at random from $Y$ with respect to the uniform probability
measure on $Y$,  then
\be
\label{discrete}
\sup_{y\in \wt Y}\|v(y)\|_V\ge  \frac{M}{8m}
\ee
with probability larger than $1-\(1-\frac{3}{4 m^{2}}\)^N$.
\end{lemma}

\noindent
{\bf Proof:}    From the definition of $m$, there exists a downward 
closed set $\Lambda$ with $\#(\Lambda)=m$ and a $V$-valued polynomial $P\in \cP_\Lambda$ such that
\be
\label{first}
\|v-P\|_{L^\infty(Y,V)}\le \frac M{4m}.
\ee
We use the Legendre polynomials  to represent $P$.   
 We denote by $(L_j)_{j\geq 0}$ the sequence of univariate Legendre polynomials 
normalized in $L^2([-1,1],\frac {dt} 2)$. Their multivariate counterparts 
\be
L_\nu(y):=\prod_{\nu_j\neq 0} L_{\nu_j}(y_j), \quad \nu\in \cF,
\ee 
are an orthonormal basis on $L^2(Y,\rho)$, where $\rho$ is the uniform probability measure 
on $Y$.  We write $P$ in its Legendre expansion
\be
\label{Legendre}
P(y)=\sum_{\nu\in\Lambda} c_\nu L_\nu(y),\quad c_\nu:=\int_Y P(y)L_\nu(y) d\rho,
\ee
where the coefficients $c_\nu$ are elements of $V$. We next invoke a result from \cite{CCMNT}
which says that for any downward closed set $\Lambda$, one has
\be
\max_{y\in Y}\sum_{\nu\in\Lambda}|L_\nu(y)|^2\le \#(\Lambda)^2.
\label{ccmnt}
\ee
Thus, it follows from the Cauchy-Schwartz inequality that for any $y\in Y$,
\be
\label{Legendre1}
\|P(y)\|_V^2\le \sum_{\nu\in\Lambda} \|c_\nu\|_V^2 \sum_{\nu\in\Lambda} | L_\nu(y)|^2\le m^2\int_Y\|P(y)\|_V^2\,  d\rho,
\ee
or equivalently
\be
\label{Legendre2}
\|P\|_{L^\infty(Y,V)}\le m\|P\|_{L^2(Y,V)}.
\ee
Now, let $S:=\{y\in Y\; :\; \|P(y)\|_V\ge \delta\}$ where $\delta:= \frac{M_1}{2m}$ with $M_1:=\|P\|_{L^\infty(Y,V)}$. Then,
\be
\label{Legendre4}
\|P\|_{L^2(Y,V)}^2\le {\anew \int_{S}\|P(y)\|_V^2\,  d\rho+\int_{S^c}\|P(y)\|_V^2\,  d\rho
\le } M_1^2\rho(S)+\delta^2.
\ee
Inserting this into \iref{Legendre2} gives
\be
\label{Legendre5}
M_1^2m^{-2}\le  M_1^2\rho(S)+\delta^2\le  M_1^2\rho(S)+m^{-2}\frac{M_1^2}4.
\ee
In other words,
\be
\label{Legendre6}
\rho(S)\ge \frac{3}{4 m^{2}}.
\ee
Suppose now that $\wt Y$ is a set formed by $N$ independent draws with respect to the uniform measure $\rho$ on $Y$.  
The probability that 
none of these draws is in
$S$ is at most  $(1-\frac{3}{4 m^{2}})^N$.  So, with probability  greater than $1- (1-\frac{3}{4 m^{2}})^N$, we have
\be
\label{legendre7}
\max_{y\in \wt Y}  \|P(y)\|_V\ge \delta =\frac{M_1}{2m} .
\ee
Accordingly, with at least the same probability, we have from \iref{first}
\be
\label{legendre7}
\max_{y\in \wt Y}  \|v(y)\|_V\ge  \delta -\frac M{4m} \ge   \frac{M_1}{2m}- \frac{M}{4m}\ge   \frac{M-\frac M{4m}}{2m}- \frac{M}{4m}\ge \frac{M}{8m},
\ee
which proves the lemma.\hfill $\Box$  
 
{\anew
\begin{remark}
\label{remmarc}
Intuitively, the above proof relies on the fact that $v$ is close to a polynomial $P$,
and that the $L^\infty$ norm of $P$ on a sufficiently fine discrete set $\t Y$ is comparable to its $L^\infty$ norm on the continuous domain $Y$. A general line of research is to look for equivalences between
discrete and continuous $L^p$ norms for a given $m$-dimensional space $X_m$
of functions, thus of the form
\be
 c_1 \|v\|_{L^p(Y)} \leq \(\frac 1 {\#(\t Y)}\sum_{y\in \t Y} |v(y)|^p\)^{1/p} \leq c_2 \|v\|_{L^p(Y)}, \quad v\in X_m.
 \label{marc}
\ee
Such results are refered to as Marcinkiewicz-type discretization theorems, see \cite{Tem} for a recent survey.
In several settings where $X_m$ consists of algebraic or trigonometric polynomials, it is known that
random sampling for $\t Y$ yields such inequalities with high probability at a sampling budget $\#(\t Y)$ that grows polynomially, and sometimes linearly, with $m$. Note that the 
inequality \iref{ccmnt} is used in \cite{CCMNT} to show that for $L^2$ norms and downward closed polynomials spaces in any dimension, the norm equivalence \iref{marc} holds with high probability for random samples of cardinality $\cO(\#(\Lambda)^2)$ up to logarithmic factors.
\end{remark}
}

\begin{remark}
\label{remcheb1}
The result in the above lemma can be improved by
sampling according the tensor product Chebychev measure, that is, with $\rho$ now defined as
\be
d\rho(y):=\bigotimes_{j\geq 1} \frac {dy_j}{\pi \sqrt{1-y_j^2}}.
\label{chebmeas}
\ee
Indeed, we can apply a similar reasoning with the $L_\nu$ replaced by the 
tensorized Chebychev polynomials $T_\nu$, for which it is proved in \cite{CCMNT} 
that 
\be
\max_{y\in Y}\sum_{\nu\in\Lambda}|L_\nu(y)|^2\le \#(\Lambda)^{2\alpha}, \quad \alpha:=\frac{\ln 3}{2\ln 2}.
\ee
for any downward closed set $\Lambda$. Therefore, the statement of the
lemma is modified as follows: if $m$ is such that $4 M_0m^{-r+\alpha}< M$ and $\wt Y$ is any finite set of $N$ points drawn at random from $Y$ with respect $\rho$,  then
\be
\label{discrete}
\sup_{y\in \wt Y}\|v(y)\|_V\ge  \frac{M}{8m^\alpha}
\ee
with probability larger than $1-\(1-\frac{3}{4 m^{2\alpha}}\)^N$.
\end{remark}
 
 \section{The main result}
 \label{sec:main}
 We are now in position to formulate our main result.  We suppose that we are given an error tolerance $\e$ and we  wish to use a greedy algorithm to construct a space $V_n$
 such that with high probability, say probability greater than $1-\eta$, we have
 \be
 \label{goal}
{ \rm dist}(\cM,V_n)\le \e,
 \ee
 with $n$ hopefully small and the off-line complexity also acceptable.  We assume that the  solution map $y\mapsto u(y)$
 belongs to $\cA^r$ for some $r>2$ and that we have an upper bound $M_0$ for $\|u\|_{\cA^r}$. 
 This assumption is known to hold in a great variety of settings of parametric PDEs, see \cite{CD}. 
 
  Given $r$ and the user prescribed $\eta$, we first define $m$ as the smallest integer such that   
 \be
  \label{defm}
 32M_0 m^{-r+2}\le \e\quad {\rm and} \quad 2^{4r+2}m^{-r}\le 1.
  \ee  
 We then define $N$  as the smallest integer such that
 \be
 \label{defN}
\(1-\frac{3}{4 m^{2}}\)^N\le \frac \eta{m^2}.  
 \ee
We consider the following greedy algorithm for finding a reduced basis.
In the first step, we make $N$ independent draws of the parameter $y$ according to
the uniform measure $\rho$. This produces a set $\wt Y_0$ of cardinality $N$.
We then use
 \be
 \label{step1}
 u^0 := u(y^0),\quad y^0:= \argmax_{y\in \wt Y_0} \|u(y)\|_V.
 \ee
 At the general step, once $u^0,\dots, u^{n-1}$ and $V_n:=\span\{u^{0},\dots,u^{n-1}\}$ have been chosen, we make $N$ independent draws of the parameter $y$ producing the set $\wt Y_{n}$ and then define
 \be
 \label{stepn}
  u^n := u(y^n),\quad y^n:= \argmax_{y\in \wt Y_n}e_n(y),
 \ee
 where $e_n(y):= \|u(y)-P_{V_n}u(y)\|_V$. In practice, each $e_n(y)$ to be computed is only approximately computed through a surrogate $\o e_{n,h}(y)$ and $u^n$ is only approximated as described in \S\ref{sec:performance}. However, we do not incorporate these facts in the analysis that follows in order to simplify the presentation.
 
 Let 
 \be
 \label{nerror}
\hat  \sigma_n:= \max_{y\in \wt Y_n} e_n(y),\quad   \sigma_n:= \max_{y\in Y} e_n(y),\quad n\ge 0,
\ee
be  the computed error and the true error for approximation of $\cM$ by $V_n$.
We terminate the algorithm at the smallest integer $n\le m^2$ for which $\hat \sigma_n\le \frac {\e}{8m}$.  
If this does not occur before $m^2$ steps we then terminate after step $n=m^2$ has been completed.   
The $V_n$ is the output of the algorithm.
 
 \begin{theorem}
 \label{maintheorem}    
 {\anew Assume that the  solution map $y\mapsto u(y)$
 belongs to $\cA^r$ for some $r>2$ with upper bound $M_0$ for $\|u\|_{\cA^r}$.
 Then,} with probability greater than $1-\eta$ the following hold for the above numerical algorithm:
 
 \noindent
 {\rm (i)} The algorithm   produces a reduced basis space $V_n$ such that $\dist(\cM,V_n)\le \e$.
 
 \noindent
 {\rm (ii)}  If  for some $s>0$ and $C_0>0$, we have    $d_n(\cM)\le C_0\max\{1,n\}^{-s}$ for all $n\ge 0$,  then the algorithm terminates in $n(\e)$ steps, where
 \be
 n(\e) \leq C\e^{-(\frac 1 s +\frac {3}{s(r-2)})},
  \label{K}
 \ee
and requires $N(\e)$ error bound evaluations, where
 \be
 N(\e)=C\e^{-\frac{2s+r+1}{s(r-2)}}(|\ln \eta|+|\ln \e|).
 \label{NK}
 \ee
The constants $C$ in the above bounds depend only on $(r,s,C_0,M_0)$.
 \end{theorem}

{\anew 
 \begin{remark}
 Note that the assumption $u\in\cA^r$ implies in particular that 
 the Kolmogorov $n$-width of $\cM$ decays at least like $n^{-r}$ and therefore
 we may assume that $s\geq r$ in the above theorem, although this is not used
 in the proof.
 \end{remark}
 }

 \noindent
 {\bf Proof:}  We first show that  with probability greater than $1-\eta$, the algorithm produces at each step 
 $k\le n$ a snapshot $u^k=u(y^k)$   
 which realizes a weak greedy algorithm, applied over {\em all} of $Y$, with parameter
 $\gamma:=\frac1{8m}$. Indeed, for any $k$,  let $v(y)= u(y)-P_{V_k}u(y)$ be the error function at the step after $V_k$ is defined.   As shown in the previous section, $v\in\cA^r$ and $\|v\|_{\cA^r}\le  M_0$.  Since the algorithm has not terminated, we have
 \be
 \label{terminate}
  \|v\|_{L^\infty(Y,V)}=\sigma_k \ge  \hat\sigma_k\ge \frac{ \e}{8m}\ge   4M_0 m^{-r+1},    
  \ee
  where the last inequality is the first condition in \iref{defm}.  Therefore,  we can apply
    Lemma \ref{slemma} to $v$ and find that 
 with probability greater than $1-\(1-\frac{3}{4 m^{2}}\)^N$, and thus from  \iref{defN} with probability greater 
 than $1-\frac{\eta}{m^2}$, we have
 \be
 \label{selection} 
\hat \sigma_k= \max_{y\in \wt Y_k}\|v(y)\|_V\ge \gamma \sup_{y\in Y} \|v(y)\|_V=\gamma \sigma_k.
 \ee
This  means that with this probability  the function  $u^{k+1}$ is a selection of the weak greedy algorithm with parameter $\gamma$.  
Since, the draws are independent  and there are at most $m^2$ sets $\wt Y_k$, the union bound
implies that with probability at least $1-\eta$, the sequence $u^1,\dots,u^n$ is a sequence that is the realization of the weak greedy algorithm with this parameter.  For the remainder of the proof, we put ourselves in the case of favorable probability.
 
 Now consider the termination of the algorithm.  If $n<m^2$, then
 \be
 \label{termination}
 \sigma_n \le  8m \hat \sigma_n\le \e,
\ee
 and so $\dist(\cM,V_n)\le \e$.    We now check the case $n=m^2$.  
 Since by assumption the solution map belongs to $\cA^r$, we know that
 \be
 d_k(\cM)\le \|u\|_{\cA^r} \max\{1,k\}^{-r}\leq M_0 \max\{1,k\}^{-r}, \quad k\geq 0.
 \ee  
 From the estimates on the performance of the weak greedy algorithm given in Theorem \ref{theo1},
 we thus know that
 \be
 \label{gperf}
 \sigma_n\le 2^{4r+1}64 m^2M_0 n^{-r}
 =2^{4r+2}32 m^{2-2r}M_0
 \le \e, 
 \ee
 where we have used the product of the two conditions in \iref{defm}.
Hence at step $n=m^2$ we have $\dist(\cM,V_n)\le \e$. Therefore,  we have completed the proof of (i).

 We next prove (ii).   So assume that  $d_n(\cM)\le C_0\max\{1,n\}^{-s}$, {\anew for some $s>0$.}
 Then, according to Theorem \ref{theo1},
 \be
 \label{wg}
 \dist(\cM,V_n)=\sigma_n\le 2^{4s+1} \gamma^{-2}C_0n^{-s}\le 2^{4s+1} 64 m^2 C_0n^{-s}, \quad n\geq 1.
 \ee
 It follows that the numerical algorithm will terminate at a $n(\e)$ with $n(\e)\le n$ where $n$ is the smallest integer that satisfies
 \be
 \label{termination}
 2^{4s+1}  64 m^2 C_0n^{-s}\le \frac{\e}{8m}.
 \ee
 Therefore,
 \be
 n(\e)\le 2^{4+\frac {10} s}C_0^{1/s} m^{\frac{3}{s}}\e^{-\frac 1 s}.
 \ee
 Using the first condition in \iref{defm}, this leads to the estimate \iref{K}
 with multiplicative constant $C:=2^{4+\frac {10} s}C_0^{1/s}(32M_0)^{\frac{3}{s(r-2)}}$.

 Finally, to execute the algorithm, we will need to draw $n(\e)$ sets $(\wt Y_0,\dots,\wt Y_{n-1})$, 
 each of them of size $N$. The total number of error bound evaluation is thus 
 \be
 N(\e)=n(\e)N.
 \ee
From the definition of $N$ in \iref{defN} we derive that 
\be
N\leq 1+\(\ln\(1-\frac{3}{4 m^{2}}\)\)^{-1}(\ln |\eta|+2 \ln \ m)\leq Cm^2(|\ln\ \eta|+\ln \ m).
\ee
Using the first condition in the definition \iref{defm} of $m$, this leads to
\be
 N\le  C \e^{-\frac{2}{(r-2)}} (|\ln \eta|+|\ln \e|). 
 \ee
 where $C$ depends on $r$ and $M_0$. Combining this with \iref{K},
 we obtain \iref{NK}, which concludes the proof of (ii).
  \hfill $\Box$
  
\begin{remark}
\label{remcheb2}
The above theorem can be improved by
sampling according the tensor product Chebychev measure \iref{chebmeas}, 
in view of Remark \ref{remcheb1}. Here, we require the solution map $y\mapsto u(y)$ 
belongs to $\cA^r$ for some $r>2\alpha:=\frac {\ln 3}{\ln 2}$. We then define $m$ as the smallest integer such that
 \be
\label{defmalt}
 32M_0 m^{-r+2\alpha}\le \e\quad {\rm and} \quad 2^{4r+2}m^{-(2\alpha-1)r}\le 1.
\ee  
and $N$ as the smallest integer such that
 \be
 \label{defNalt}
\(1-\frac{3}{4 m^{2\alpha}}\)^N\le \frac \eta{m^{2\alpha}}.  
 \ee
We terminate the algorithm at the smallest integer $n\le m^{2\alpha}$ for which $\hat \sigma_n\le \frac {\e}{8m^\alpha}$.  
With the exact same proof, we reach the statement as Theorem \ref{maintheorem}, however with a number of step
 \be
 n(\e) \leq C\e^{-(\frac 1 s+\frac {3\alpha}{s(r-2\alpha)})},
  \label{K2}
 \ee
and a number of error bound evaluations
 \be
 N(\e)=C\e^{-\frac{2s\alpha+r+\alpha}{s(r-2\alpha)}}(|\ln \eta|+|\ln \e|).
 \label{NK2}
 \ee
\end{remark}
 
 Let us comment on the difference in performance between the above algorithm using random sampling 
and the greedy algorithm based on using an $\e$-net
for the solution manifold as a training set. We aim at a target accuracy $\e$, and
assume that the $n$-widths of the solution manifold decay like $d_n(\cM)\le Cn^{-s}$.

Then, the approach based on an $\e$-net
constructs a reduced basis space $V_n$ of optimal dimension
\be
n=n(\e)\sim \e^{-1/s},
\ee
and the total number of error bound evaluations is at best of the order
\be
N(\e)\sim n(\e) \e^{-1/s},
\ee
each of them having a cost ${\rm Poly}(n)={\rm Poly}(\e^{-1})$, resulting in a prohibitive offline cost 
${\rm Poly}(\e^{-1})\e^{-1/s}$. Note that one should also add the cost of evaluating the 
$n(\e)$ reduced basis functions on the finite element space of large dimension $n_h>>1$.

In contrast, the approach based on random sampling constructs a 
reduced basis space $V_n$ of sub-optimal dimension
\be
n(\e) \leq C\e^{-(\frac 1 s +\frac {3}{s(r-2)})},
\label{neps}
\ee
but the total number of error bound evaluation is now of the order
\be
N(\e)\sim \e^{-\frac{2s+r+1}{s(r-2)}}(|\ln \eta|+|\ln \e|),
\label{Neps}
\ee
where $\eta$ is the probability of failure.

In summary, while our approach allows for a dramatic reduction in the offline
cost, it comes with a loss of optimality in the performance of reduced basis spaces
since $n(\e)$ scales with $\e^{-1}$ with an exponent larger than $\frac 1 s$. 
In particular, this affects the resulting online cost.
Inspection of the proof of the main theorem reveals that this
loss comes from the fact that the greedy selection from the random set can 
only be identified to a weak-greedy algorithm with a parameter $\gamma=\frac 1 {8m}$
which instead of being fixed becomes small as $m$ grows, or equivalently as $\e$ decreases. 
Let us still observe that the above perturbation
of $\frac 1 s$ by $\frac {3}{s(r-2)}$ becomes neglectible as $r$ gets larger.
We have also see that this perturbation can be reduced to $\frac {3\alpha}{s(r-2\alpha)}$
by sampling randomly according to the Chebychev measure.

This leaves open the question of finding a sampling strategy for the training set which 
lead to reduced basis of optimal complexity $n(\e)\sim \e^{-1/s}$ and where the
number of error bound evaluation in the offline stage remains polynomial in $\e^{-1}$.

{\anew
\section{Numerical illustration}
\label{sec:test}

The results that we have obtained in the previous section can be rephrased
in the following way: a polynomial rate of decay of the error achieved by the greedy algorithm
in terms of reduced basis space dimension $n$ can be maintained 
when using a random training set of cardinality $N$ that scales polynomially in $n$.
More precisely, in view of \iref{neps} and \iref{Neps}, a sufficient scaling is
\be
N\sim n^{\beta^*}, \quad \beta^*=\frac{2s+r+1}{s(r-2)} \(\frac 1 s +\frac {3}{s(r-2)}\)^{-1}
\ee
independently of the parametric dimension $d$. Note that we obviously
have 
$$
\beta^*\geq \frac {2s+r+1}{3} \geq \frac{3+2s}{3}>1.
$$

In this section, we illustrate these findings through the following numerical test: we consider
the elliptic diffusion equation \iref{ellip} set on the square domain $D=]0,1[^2$, and with
affine parametrization \iref{affine} where $\bar a=1+\delta$ for some $\delta>0$ and the $\psi_j$ are given by
\be
\psi_j=a_j \chi_{D_j}
\ee
where $\{D_1,\dots,D_d\}$ is a uniform partition of $D$ into $d=k\times k$ 
squares of equal size. We take $k=8$ and therefore $d=64$, and take 
\be
a_j=j^{-t},
\ee
for some $t>0$. The effect of taking $t$ larger is to raise the
anisotropy of the parametric dependence. The results from \cite{BCM} show that this
is directly reflected by a rate $n^{-r}$ of sparse polynomial approximation of the
parameter to solution map, for any $r<t-\frac 1 2$,
and therefore on the rate of decay of the $n$-width $d_n(\cM)$. The effect of taking $\delta$ closer
to $0$ is to make the problem more degenerate as $y_1$ gets close to $-1$.
 
 We test the performance $\sigma_n$ of the reduced basis spaces generated
 by the greedy algorithm using random training sets $\wt\cM$ of cardinality
 \be
 N=N(n)=\lfloor n^{\beta}\rfloor,
 \ee
 for some $\beta\geq 1$. Note that the case $\beta=1$ amounts in selecting the $n$
 reduced basis element completely at random since all the elements from $\wt \cM$
 are necessarily chosen. We expect the performance to improve as $\beta$ becomes
 larger since we then perform a particular selection of the $n$ reduced basis elements 
 within $\cM$ through the greedy algorithm. On the other hand, our theoretical results indicate
 that a fixed value $\beta=\beta^*$ should be sufficient to ensure that the algorithm performs
almost as good as if the selection process took place on the whole of $\cM$. 

In our numerical test, we took $\delta=10^{-2}$, that is $\bar a=1.01$.
As to the parametric dimension, we considered $d=16$ and $d=64$
that correspond to subdivisions of $D$ into  $4\times 4$ and $8\times 8$
subdomains, respectively. As to the decay of $a_j$,
we test the values $t=1$ and $t=2$. Finally 
for the growth of the training sample size $N(n)$, we test the values
\be
\beta=1,1.25,1.5,1.75,2.
\ee
The error curves of the reduced basis approximation, averaged over $20$ realizations of the
random training sample, for the various choices of $d$, $t$ and $\beta$,
are displayed on Figures \ref{fig16} and \ref{fig64}. 

\begin{figure}
\begin{tabular}{cc}
	\includegraphics[width=8cm]{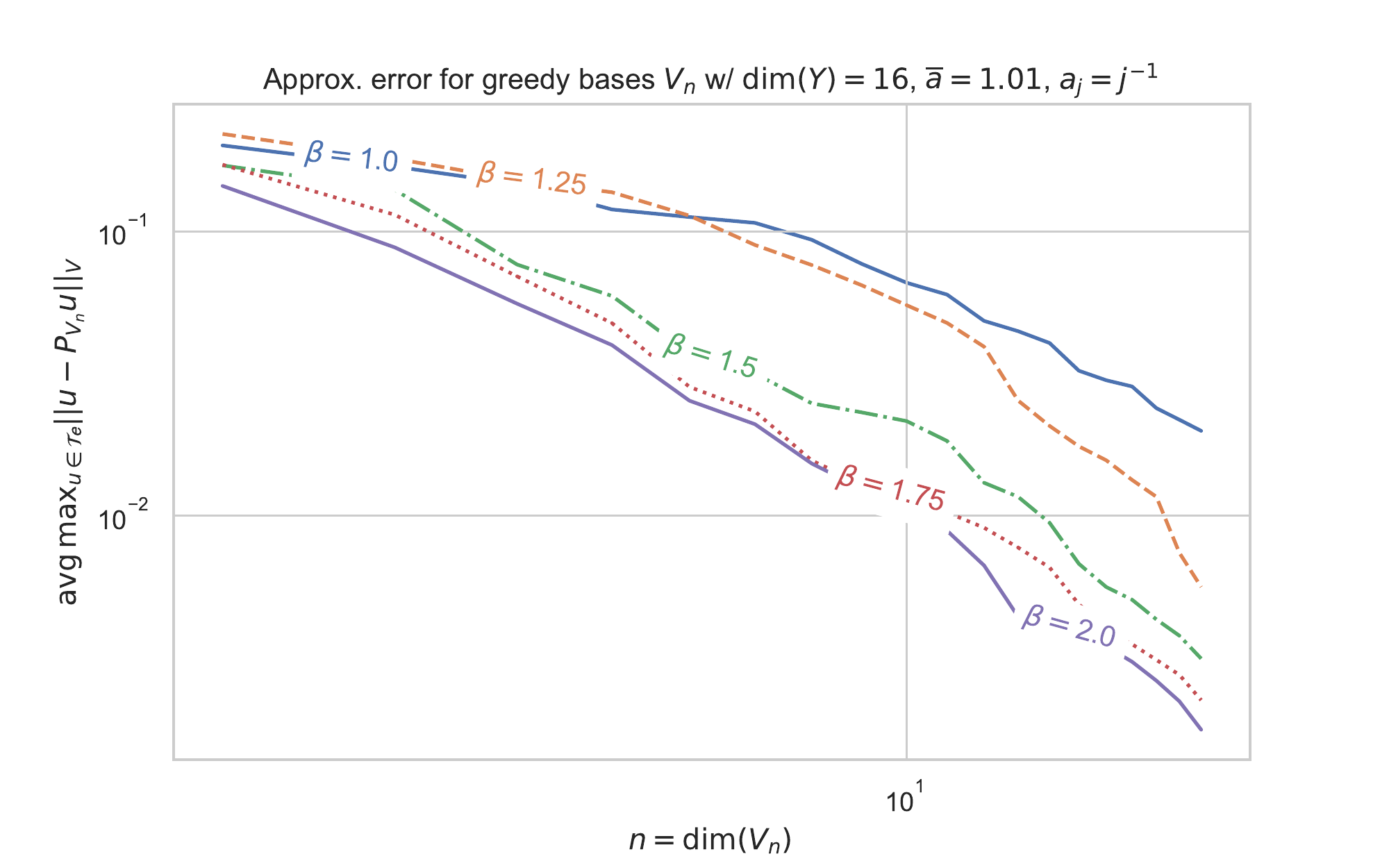} & \includegraphics[width=8cm]{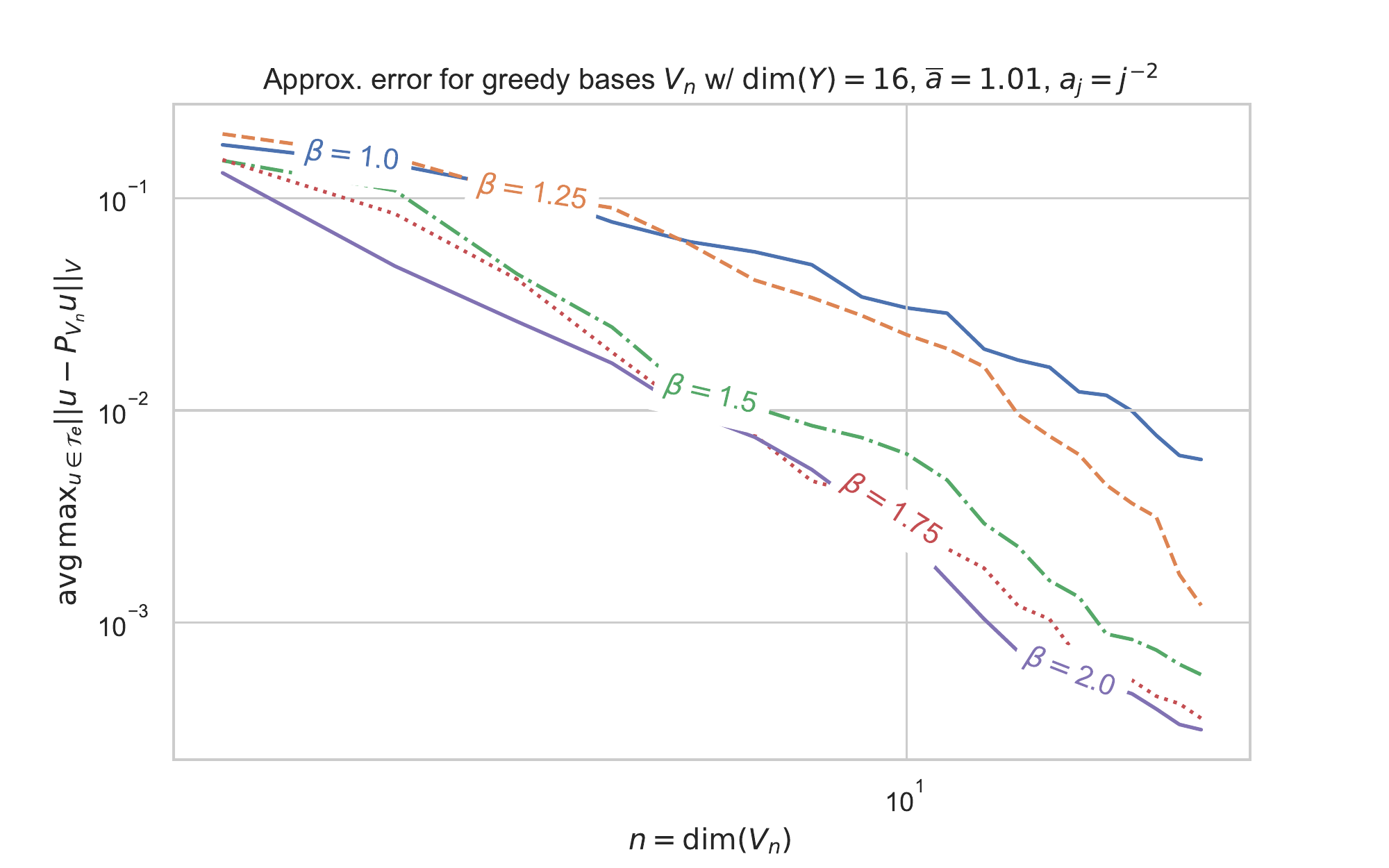} \\
	\end{tabular}
	\caption{Reduced basis approximation with $d=16$ and decay rate $t=1$ or $t=2$.}\label{fig16}
\end{figure}

\begin{figure}
\begin{tabular}{cc}
	\includegraphics[width=8cm]{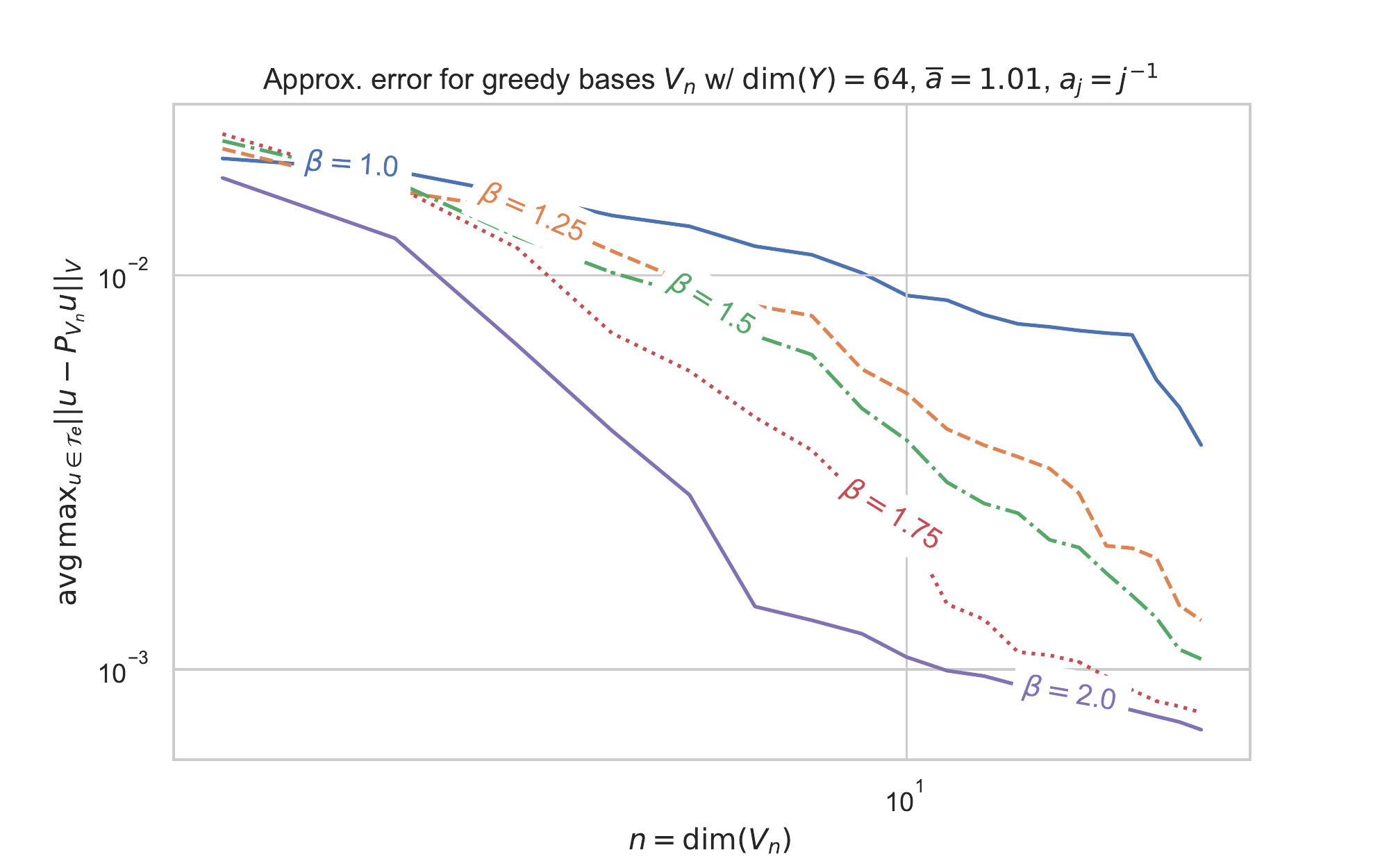} & \includegraphics[width=8cm]{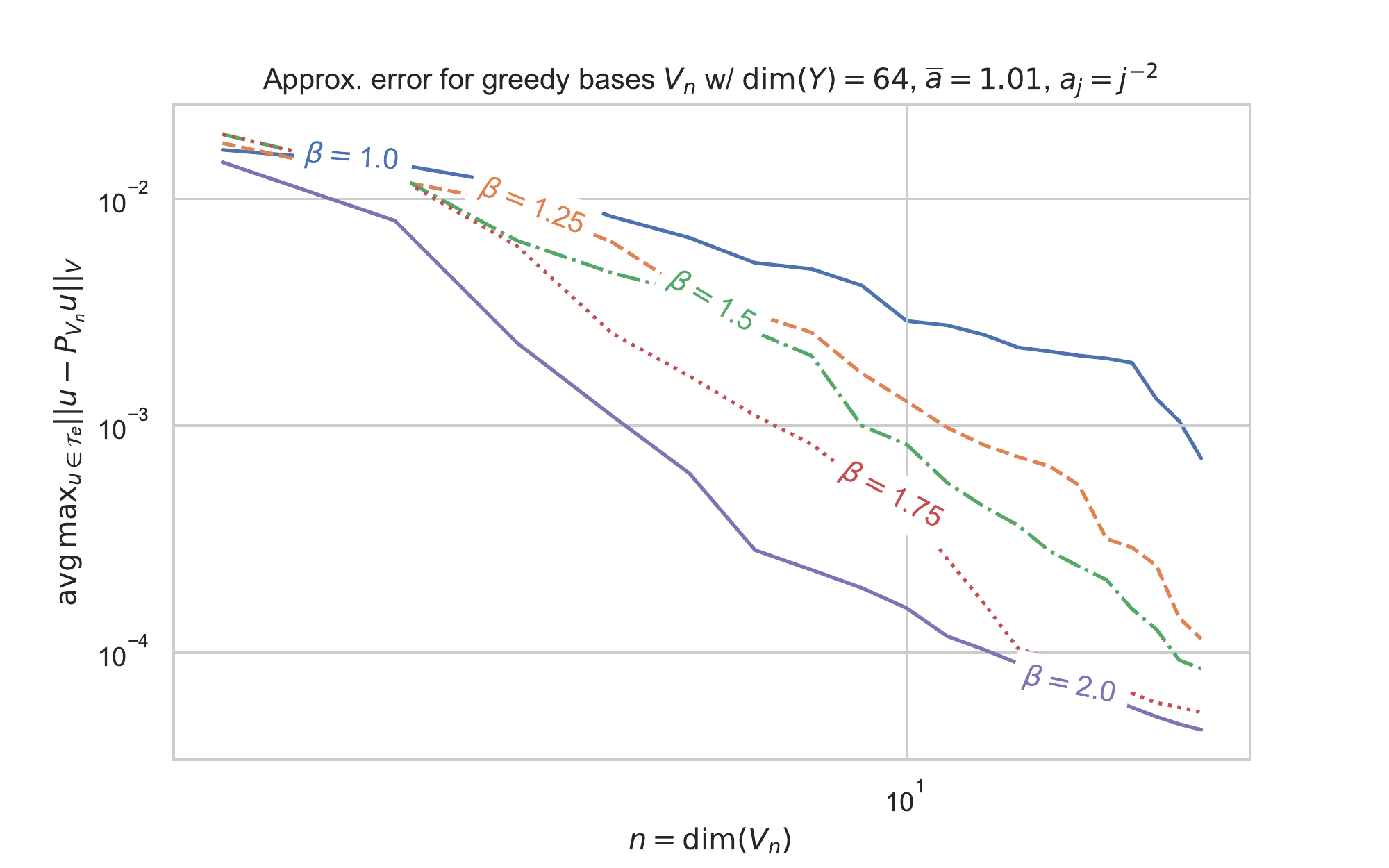} \\
	\end{tabular}
	\caption{Reduced basis approximation with $d=64$ and decay rate $t=1$ or $t=2$.}\label{fig64}
\end{figure}

 As predicted by the theoretical approximation results, the convergence rate of the reduced basis
 method improves as $t$ gets larger. The observed convergence rates (for the highest
 value of $\beta$) are closer to the value $t$ than to $t-1/2$ which can be rigourously
 established from the polynomial approximation rates. This reflects the fact that the 
 reduced basis approximations perform generally better than 
 sparse polynomial approximations.

We also note that, for the same value of $t$, the errors are smaller in the higher parametric dimension $d=64$
than for $d=16$. This apparent paradox can be explained:
in both cases, the most active variables are the first ones ($y_1$, then $y_2$,..)
yet they are associated to domains $D_j$ of smaller size in the high dimensional case,
therefore having less impact on the variation of the solution with these variables.

As expected, we observe that the error curve behaves better as we increase the value
of $\beta$ but we observed that this phenomenon stagnates at $\beta=2$.
This hints that the scaling $N(n)=n^2$ is practically sufficient 
to ensure in this case the optimal convergence behaviour which
would be met with a very rich training set, for example an $\e$-net. The value $\beta=2$ is
much smaller than the value $\beta^*$ given by the theoretical analysis which is thus
too much pessimistic.
}

\vskip .2in

\noindent 
Albert Cohen,  Laboratoire Jacques-Louis Lions, Sorbonne Universit\'e, Paris, France.
\vskip .1in
\noindent
Wolfgang Dahmen, Department of Mathematics,  University of South Carolina,  Columbia, SC 29208
\vskip .1in
\noindent
Ronald DeVore,  Department of Mathematics, Texas A$\&$M University,  College Station, TX  77840
\vskip .1in
\noindent 
James Nichols,  Laboratoire Jacques-Louis Lions, Sorbonne Universit\'e, Paris, France.
\end{document}